\documentclass{crm-eca}

\usepackage{amssymb}

\usepackage{graphicx}

\usepackage{HamSysprepr}  

\newcommand\R{\mathbb R}
\newcommand\T{\mathbb T}
\newcommand\Z{\mathbb Z}

\newcommand\ee{\mathrm e}

\newcommand\eps\varepsilon

\newcommand\J{\mathcal J}
\newcommand\Lc{\mathcal L}
\newcommand\M{\mathcal M}

\newcommand\W{\mathcal W}
\newcommand\Tc{\mathcal T}

\newcommand\p[1]{\left(#1\right)}
\newcommand\pq[1]{\left[#1\right]}
\newcommand\pp[1]{\left\{#1\right\}}
\newcommand\scprod[2]{\left\langle#1,#2\right\rangle}
\newcommand\abs[1]{\left|#1\right|}

\newcommand\tl{\tilde}

\newcommand\df{{\rm d}}
\newcommand\Df{{\rm D}}
\newcommand\Ord{\mathcal O}
\newcommand\rint{\mathop{\rm rint}}

\newcommand\mmatrix[4]{\p{\begin{array}{cc}#1&#2\\[3pt]#3&#4\end{array}}}
\newcommand\symmatrix[3]{\mmatrix{#1}{#2}{#2}{#3}}

\newcommand\ut{{\rm u}}  
\newcommand\st{{\rm s}}  

\newtheorem{theorem}{Theorem}

\theoremstyle{definition}

\theoremstyle{remark}

\begin{document}

\papertitle[Exponentially small splitting with quadratic frequencies]{A methodology for obtaining asymptotic estimates for the exponentially
       small splitting of separatrices to whiskered tori with quadratic frequencies
}


\paperauthor{Amadeu Delshams}
\paperaddress{Departament de Matem\`atica Aplicada I, 
    Universitat Polit\`ecnica de Catalunya,
    } 
\paperemail{amadeu.delshams@upc.edu, pere.gutierrez@upc.edu}

\paperauthor{Marina Gonchenko}
\paperaddress{Institut f\"ur Mathematik, Technische Universit\"at Berlin,
    } 
\paperemail{gonchenk@math.tu-berlin.de}

\paperauthor{Pere Guti\'errez}


\paperthanks{This work has been partially supported by the Spanish
     MINECO-FEDER grant MTM2012-31714 and the Catalan grant 2009SGR859.
     The author~MG has also been supported by the
     DFG~Collaborative Research Center TRR~109
     ``Discretization in Geometry and Dynamics''.}

\makepapertitle




\subsection*{Introduction}

The aim of this work is to provide asymptotic estimates for the splitting of
separatrices in a perturbed 3-degree-of-freedom Hamiltonian system, associated
to a 2-dimensional whiskered torus (invariant hyperbolic torus) whose frequency
ratio is a quadratic irrational number. We show that the dependence of the
asymptotic estimates on the perturbation parameter is described by
some functions which satisfy a periodicity property, and whose
behavior depends strongly on the arithmetic properties of the frequencies.

First, we describe the Hamiltonian system to be studied.
It is also considered in \cite{DelshamsGS04},
as a generalization of the famous Arnold's example \cite{Arnold64},
and provides a model for the behavior of a nearly-integrable
Hamiltonian system in the vicinity of a single resonance
(see \cite{DelshamsG01} for a motivation).
In canonical coordinates
$(x,y,\varphi,I)\in\T\times\R\times\T^2\times\R^2$,
we consider a perturbed Hamiltonian
\begin{eqnarray}
  \label{eq:HamiltH}
  &&H(x,y, \varphi, I) = H_0 (x,y, I) + \mu H_1(x, \varphi).
\\
  \label{eq:HamiltH0}
  &&H_0 (x, y, I) =
  \langle \omega_\eps, I\rangle + \frac{1}{2} \langle\Lambda I, I\rangle
  +\frac{y^2}{2} + \cos x -1,
\\
  \label{eq:HamiltH1}
  &&H_1 (x, \varphi)
  =\cos x\cdot\sum_{k_2\ge0}
   \ee^{-\rho |k|}\cos(\langle k,\varphi\rangle-\sigma_k).
\end{eqnarray}
For the integrable Hamiltonian $H_0$, we consider
a vector of \emph{fast frequencies}
\begin{equation}
  \omega_\eps = \frac\omega{\sqrt{\eps}}\;,
  \qquad
  \omega=(1,\Omega), 
\label{eq:omega_eps}
\end{equation}
where the frequency ratio $\Omega$ is a \emph{quadratic} irrational number.
In this way, our system has two parameters $\eps>0$ and~$\mu$,
but we assume them linked by a relation of kind $\mu=\eps^p$, $p>0$
(the smaller $p$ the better).
Thus, if we consider $\eps$ as the unique parameter, we have
a \emph{singular} or \emph{weakly hyperbolic} problem for $\eps\to 0$
(see \cite{DelshamsG01} for a discussion about singular and regular problems).

On the other hand, notice that $H_0$ consists of a classical pendulum 
and $2$~rotors (in the coordinates $x,y$ and $\varphi,I$ respectively).
Then, we see that $H_0$ has a family of 2-dimensional whiskered tori,
with coincident whiskers (invariant manifolds).
Such tori can be indexed by the (constant) action $I$,
and have frequency vectors $\omega_\eps+\Lambda I$.
We assume that the matrix $\Lambda$ is such that
the condition of \emph{isoenergetic nondegeneracy} is satisfied
(see for instance \cite{DelshamsGS04}).
Among the tori, we fix our attention on the torus given by $I=0$,
\[
  \Tc_0: \qquad (0,0,\theta,0), \quad \theta\in\T^2,
\]
whose inner flow is given, in this parameterization,
by $\dot\theta=\omega_\eps$.
This torus has a homoclinic whisker
(i.e.~coincident stable and unstable whiskers),
\[
  \W_0: \qquad (x_0(s),y_0(s),\theta,0), \quad s\in\R,\ \theta\in\T^2,
\]
where $x_0(s)=4\arctan\ee^s$, $y_0(s)=2/\cosh s$
(the upper separatrix of the classical pendulum).
The inner flow on $\W_0$ is given by $\dot s=1$, $\dot\theta=\omega_\eps$\,.

Concerning the perturbation $H_1$, it is given by a constant $\rho>0$
(the complex width of analyticity in the angles~$\varphi$),
and phases $\sigma_k$ that, for the purpose of this work, can be chosen
arbitrarily.

Under the hypotheses described,
the \emph{hyperbolic KAM theorem} (see for instance \cite{Niederman00})
can be applied to the perturbed
Hamiltonian~(\ref{eq:HamiltH}--\ref{eq:HamiltH1}).
We have that, for $\mu\ne0$ small enough, the whiskered torus $\Tc_0$ persists
with some shift and deformation giving rise to a perturbed torus $\Tc$,
with perturbed local stable and unstable whiskers.

Such local whiskers can be extended to global whiskers $\W^\st$, $\W^\ut$
but, in general, for $\mu\ne0$ they do not coincide anymore,
and one can introduce a \emph{splitting function}
giving the distance between the whiskers
in the directions of the actions $I\in\R^2$:
denoting $\J^{\st,\ut}(\theta)$ parameterizations
of a transverse section of both whiskers, one can define
\ $\M(\theta):=\J^\ut(\theta)-\J^\st(\theta)$, \ $\theta\in\T^2$.
In fact, this function turns out to be the gradient
of the (scalar) \emph{splitting potential}:
\ $\M(\theta)=\nabla\Lc(\theta)$
\ (see \cite[\S5.2]{DelshamsG00}, and also~\cite{Eliasson94}).

In~(\ref{eq:omega_eps}), we deal  with the following  24~quadratic numbers 
\begin{equation}
  [\bar{1}]\,, [\bar{2}]\,, \;\ldots, \; [\overline{{13}}]\,, [\overline{{1,2}}]\,,\;\ldots, \; [\overline{{1,12}}]\,,
\label{eq:24numbers}
\end{equation}
where we denote a quadratic number according to its periodic part
in the continued fraction (see~(\ref{eq:cont_frac})).

Next, we establish the \emph{main result} of this work,
providing two types of \emph{asymptotic estimates} for the splitting,
as $\eps\to0$.
One one hand, we give an estimate for the \emph{maximal splitting distance},
i.e.~for the maximum of $\abs{\M(\theta)}$, $\theta\in\T^2$.
On the other hand, we show that for most values of $\eps\to0$
there exist 4~transverse homoclinic orbits, associated to
simple zeros $\theta_*$ of $\M(\theta)$
(i.e.~nondegenerate critical points of $\Lc(\theta)$)
and, for such homoclinic orbits, we obtain an estimate
for the \emph{transversality} of the splitting, given by the
minimum eigenvalue (in modulus) of the matrices $\Df\M(\theta_*)$.

We use the notation $f \sim g$
if we can bound $c_1 |g| \leq |f| \leq c_2 |g|$ with positive constants
$c_1, c_2$ not depending on $\eps$, $\mu$.

\begin{theorem}\label{thm:main}
Assume the conditions described above for the
Hamiltonian~\mbox{(\ref{eq:HamiltH}--\ref{eq:HamiltH1})},
and that $\eps$ is small enough and $\mu=\eps^p$, $p>3$.
Then, there exist continued functions $h_1(\eps)$ and $h_2(\eps)$
(defined in~(\ref{eq:h12})),
periodic in $\ln\eps$ and satisfying $1\le h_1(\eps)\le h_2(\eps)$,
and a positive constant $C_0$ (given in~(\ref{eq:C0})), such that:
\begin{itemize}
\item[\rm(a)]
for the maximal splitting distance, we have the estimate
$$
\max_{\theta\in \T^2} |\M(\theta)|
\sim \frac{\mu}{\sqrt{\eps}}
\exp \left\{- \frac{C_0 h_1 (\eps)}{\eps^{1/4}}\right\};
$$
\item[\rm(b)]
the splitting function $\M(\theta)$ has exactly 4 zeros $\theta_*$,
all simple, for all $\eps$ except for a small neighborhood
of a finite number of geometric sequences of $\eps$;
\item[\rm(c)]
at each zero $\theta_*$ of $\M(\theta)$,
the minimal eigenvalue of $\Df\M(\theta_*)$ satisfies
\begin{equation*}\label{eq:minimaleig}
  m_* \sim \mu \eps^{1/4}
  \exp \left\{- \frac{C_0 h_2 (\eps)}{\eps^{1/4}}\right\}.
\end{equation*}
\end{itemize}
\end{theorem}

For the proof of this theorem, we apply the \emph{Poincar\'e--Melnikov method},
which provides a first order approximation
\begin{equation}\label{eq:melniapprox}
  \M(\theta)=\mu\nabla L(\theta)+\Ord(\mu^2)
\end{equation}
in terms of the \emph{Melnikov potential}, which can be defined by integrating
the perturbation $H_1$ along the
trajectories of the unperturbed homoclinic whisker $\W_0$\,:
\begin{equation}\label{eq:L}
  L(\theta):= - \int_{-\infty}^{\infty}
  H_1(x_0(t),\theta+\omega_\varepsilon t)\,\df t.
\end{equation}
Since this first order approximation is exponentially small in $\eps$,
in principle the approximation~(\ref{eq:melniapprox}) cannot be directly
applied in our singular problem with $\mu=\eps^p$.
However, using suitable bounds for the error term $\Ord(\mu^2)$,
given in \cite{DelshamsGS04}, one can see that for $p>3$ the
first order approximation given by the Melnikov potential overcomes
the error term and provides the right asymptotic estimates for the splitting.
Such estimates come from the size of \emph{dominant harmonics}
in the Fourier expansion of~(\ref{eq:L}),
and studying their dependence on $\eps$.
More precisely, to estimate the maximal splitting
one dominant harmonic is enough
and, to estimate the transversality of the splitting,
two dominant harmonics are required
(excluding the values of $\eps$ such that the second and third
harmonics are of the same magnitude, which could give rise to bifurcations
in the homoclinic orbits and would require a further study).

The remaining sections of this work are devoted to the definition of the
functions $h_1(\eps)$ and $h_2(\eps)$, putting emphasis on their
dependence on the arithmetic properties of the quadratic number $\Omega$.

\subsection*{Continued fractions and resonant sequences}

We review briefly the technique developed in \cite{DelshamsG03}
for studying the resonances of quadratic frequencies.
Let $0<\Omega<1$ be a quadratic irrational number.
It is well-known that it has an infinite continued fraction
\begin{equation}\label{eq:cont_frac}
  \Omega=[a_1,a_2,a_3,\ldots]
  =\frac{1}{a_1+\dfrac{1}{a_2+\dfrac{1}{a_3+\cdots}}}\,,
  \qquad
  a_n\in\Z^+,\ n\ge1
  \quad \mbox{(and $a_0=0$)},
\end{equation}
that is that is \emph{eventually periodic}, i.e., periodic starting at some $a_l$. 
For a purely $m$-periodic
continued fraction $\Omega=[\overline{a_1,\ldots,a_m}]$ we introduce the matrix 
$$
U =(-1)^m A_1^{-1}\cdots A_m^{-1}, \textrm{ where } A_l =\symmatrix{a_l}10, \;\;\; l=1,\ldots, m. 
$$

It is well-known that quadratic vectors satisfy
a \emph{Diophantine condition}
\begin{equation*}
|\langle k, \omega\rangle| \geq \frac{\gamma}{|k|},
\;\;\; \forall k\in \Z^2\setminus\{0\}.
\label{eq:DiophCond}
\end{equation*}
With this in mind, 
we define the \emph{``numerators''}
\begin{equation}\label{eq:numerators}
  \gamma_k := |\langle k, \omega\rangle|\cdot|k|,
  \qquad
  k\in\Z^2\setminus\pp0
\end{equation}
(for integer vectors, we use the norm $\abs\cdot=\abs\cdot_1$).
Our aim is to find the integer vectors $k$ which give
the smallest values $\gamma_k$,
we call such vectors $k$ the \emph{primary resonances}.

All vectors $k\in\Z^2\setminus\pp0$
with $\abs{\scprod k\omega}<1/2$
are subdivided into \emph{resonant sequences}:
\begin{equation}
s(j,n) := U^n k^0(j),
\qquad n=0,1,2,\ldots
\label{eq:sjn}
\end{equation}
where the initial vector $k^0(j)=(-\rint(j\Omega),j)$, $j\in \Z^+$, satisfies
\begin{equation}\label{eq:primit}
 \frac{1}{2\lambda} < |\langle k^0(j), \omega\rangle| < \frac{1}{2},
\end{equation}
$\lambda$ being the eigenvalue of $U$ with $\lambda>1$. 
For each $j\in \Z^+$ satisfying~(\ref{eq:primit}),
it was proved in \cite[Th.~2]{DelshamsG03} (see also \cite{DelshamsGG13})
that, asymptotically, the resonant sequence $s(j,n)$
exhibits a geometric growth 
and the sequence $\gamma_{s(j,n)}$ has a limit $\gamma^*_j$:
\begin{equation}\label{eq:gammajast}
|s(j,n)| = K_j \lambda^n + \Ord(\lambda^{-n}), \;\;\; \gamma_{s(j,n)} = \gamma_j^* + \Ord(\lambda^{-2n}), \textrm{ as } n\to \infty,  
\end{equation}
where $K_j$ and $\gamma_j^*$ can be determined explicitly for each resonant sequence (see explicit formulas in 
\cite{DelshamsG03}). 
We select the minimal of $\gamma_j^*$:
\begin{equation}
\gamma^*:=\liminf_{|k|\to \infty} \gamma_k
= \min_{j}\gamma^*_j = \gamma^*_{j_0}>0.
\label{eq:quad_gamj0}
\end{equation}
The integer vectors of the corresponding sequence $s(j_0, n)$
are \emph{the primary resonances},
and we call the \emph{secondary resonances}  the integer vectors
belonging to any of the remaining {resonant} sequences $s(j,n)$, $j\neq j_0$. We also call by 
\emph{the main secondary resonances} the sequence $s(j_1,n)$ which is linearly independent with $s(j_0,n)$ and gives the smallest limit 
$\gamma^*_{j_1}$ among the secondary resonances.

\subsection*{The functions $h_1(\eps)$ and $h_2(\eps)$}
Taking into account the form of $H_1$ in~(\ref{eq:HamiltH1}), we present the Melnikov potential (\ref{eq:L}) in its Fourier expansion. 
Using~(\ref{eq:omega_eps}) and~(\ref{eq:numerators}),
we present the coefficients in the form
\begin{equation}
\label{eq:alphabeta}
L_k = \frac{2\pi |\langle k, \omega_\varepsilon\rangle|
\,\ee^{-\rho |k|}}{\sinh |\frac{\pi}{2}
\langle k, \omega_\varepsilon\rangle|}
    = \alpha_k\,\ee^{- \beta_k},
\qquad
\alpha_k \approx \frac{4 \pi\gamma_k}{|k|\sqrt{\varepsilon}}\,,
\quad
\beta_k =\rho |k| + \frac{\pi \gamma_k}{2 |k|\sqrt{\varepsilon}}\,.
\end{equation}
For any given $\eps$, we find the dominant harmonics $L_k(\eps)$ which
correspond essentially to the smallest exponents
$\beta_k(\eps)$. 

The exponents $\beta_k(\eps)$ in~(\ref{eq:alphabeta}) can be
presented in the form
\begin{equation}\label{eq:gk_quad}
  \beta_k(\eps) = \frac{C_0}{\eps^{1/4}}\,g_k (\eps),
  \qquad
  g_k (\eps):=
  \frac{\tl\gamma_k^{1/2}}2
  \pq{\p{\frac{\eps}{\eps_k}}^{1/4}+\p{\frac{\eps_k}{\eps}}^{1/4}},
  \end{equation}
where
\begin{equation}\label{eq:C0}
  \eps_k:=D_0\,\frac{\tl\gamma_k^{\,2}}{\abs k^4}\,,
  \qquad
  \tl\gamma_k=\frac{\gamma_k}{\gamma^*},
  \qquad
C_0=(2\pi\rho\gamma^*)^{1/2},
  \qquad
  D_0=\p{\frac{\pi\gamma^*}{2\rho}}^2.
\end{equation}
Since the coefficients $L_k$ are exponentially small in $\eps$,
it is more convenient to work with the functions $g_k$,
whose smallest values correspond  to the largest $L_k$.
To this aim, it is useful to consider the graphs of the functions
$g_k(\eps)$, $k\in\Z^2\setminus\{0\}$,
in order to detect the minimum of them for a given value of $\eps$.
\begin{figure}[b!]
\hspace*{-0.9cm}
  \includegraphics[width=1.05\textwidth, height=0.3\textwidth]{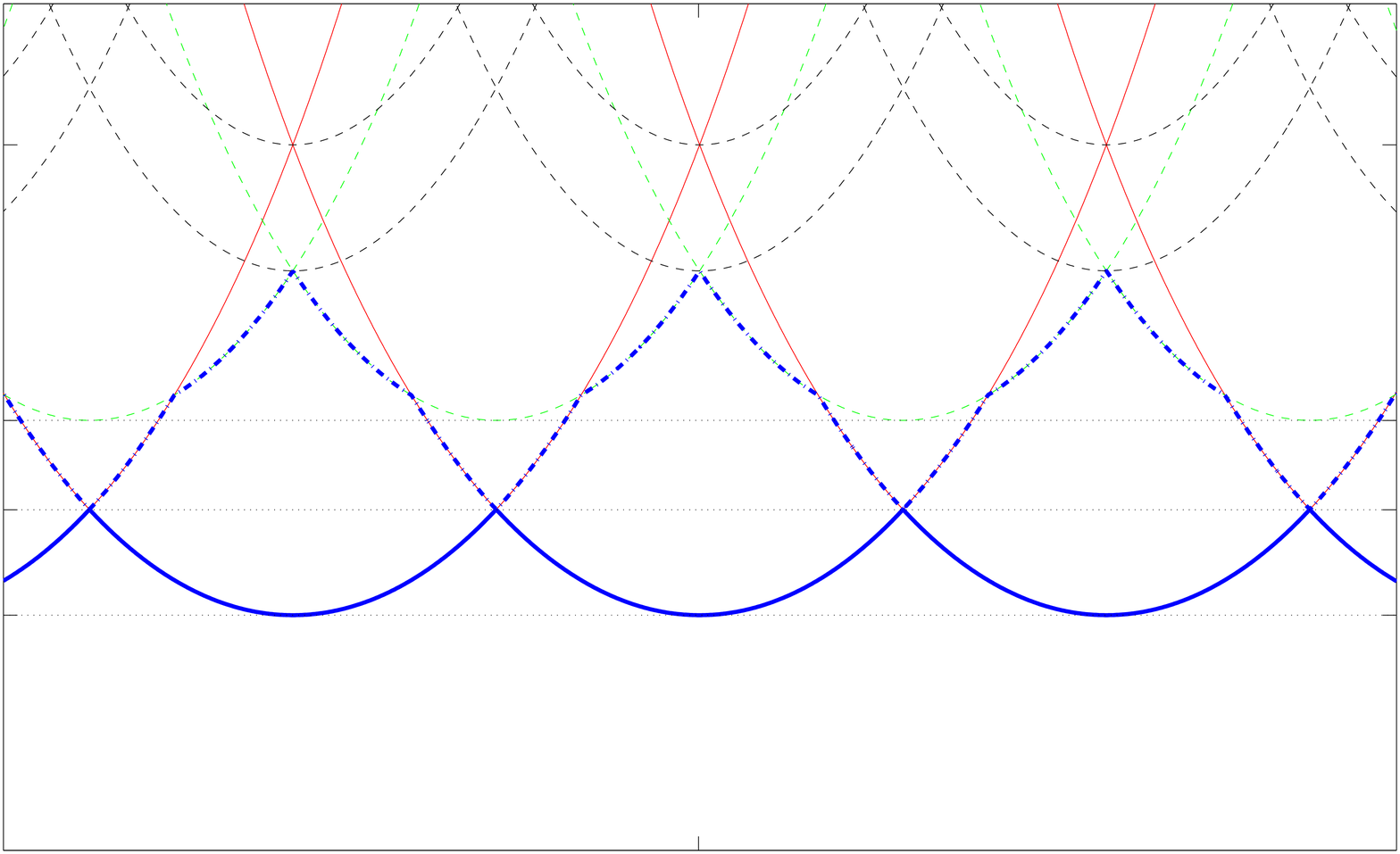}
  \caption{{\small\emph{Graphs of the functions $h_1(\eps)$ (solid blue) and $h_2(\eps)$ (dash-dot blue) for $[\overline{1,2}]=\sqrt{3}-1$. 
    The red lines are the primary functions $g_{s(j_0,n)}(\eps)$, and the green lines correspond to the main secondary functions $g_{s(j_1,n)}(\eps)$.}}}
\label{fg:qn[1p2]}
\end{figure}

We know from~(\ref{eq:gk_quad}) that the functions $g_k(\eps)$
have their minimum at $\eps=\eps_k$ and the corresponding minimal values are
$g_k(\eps_k) = \tilde{\gamma}_k^{1/2}$.
For the integer vectors $k=s(j,n)$ belonging to a resonant sequence~(\ref{eq:sjn}),
using the approximations~(\ref{eq:gammajast}), we have
$$
\eps_{s(j,n)} \approx \frac{D_0(\tl\gamma^*_j)^2}{K_j^{\,4}\lambda^{4n}}, \;\;\; 
g_{s(j,n)}(\eps) \approx \frac{(\tl\gamma^*_j)^{1/2}}2
  \pq{\p{\frac{\eps}{\eps_{s(j,n)}} }^{1/4}+\p{\frac{\eps_{s(j,n)}}{\eps}^{1/4}}},  \textrm{ as } n\to\infty.
$$
Taking into account such approximations, we have a periodic behavior of the functions with respect to $\ln\eps$, as we see in Figure~\ref{fg:qn[1p2]}
(where a logarithmic scale for $\eps$ is used).

We define, for any given $\eps$, the function $h_1(\eps)$ and $h_2(\eps)$ 
as 
\begin{equation}\label{eq:h12}
 \displaystyle h_1(\eps):=\min_k g_k(\eps)=g_{S_1}(\eps),
\;\;\;  
\displaystyle    h_2 (\eps):=
  \min\limits_{k\ {\rm lin. indep. of}\;S_1}g_k(\eps)
= g_{S_2} (\eps),
  \end{equation}
with some integer vectors $S_1(\eps)$ and $S_2(\eps)$ realizing such minima.
The functions are continuous and $4 \ln \lambda$-periodic in $\ln \eps$. It turns out that for the 24 quadratic numbers (\ref{eq:24numbers}), 
the integer vector $S_1(\eps)$ providing $h_1(\eps)$ always corresponds to a primary resonance, defined in~(\ref{eq:quad_gamj0}). 
On the other hand, the vector $S_2(\eps)$ providing $h_2(\eps)$ may correspond to primary or main secondary resonances in different intervals of $\eps$
(see~Figure~\ref{fg:qn[1p2]} for an illustration for the number $[\overline{1,2}]=\sqrt{3}-1$). 
There is a finite number of geometric sequences of $\eps$, where a change 
in $S_2(\eps)$ occurs. These points require a special study for the transversality and they are excluded in Theorem~\ref{thm:main}.









\end{document}